\documentclass{amsart}
\usepackage{algorithmic}
\usepackage{algorithm}
\usepackage{hyperref}
\usepackage[demo]{graphicx}
\usepackage{array}
\usepackage{tabularx}

\newtheorem{theorem}{Theorem}[section]
\newtheorem{lemma}[theorem]{Lemma}

\theoremstyle{definition}

\newtheorem{example}[theorem]{Example}

\usepackage{xcolor}
\theoremstyle{remark}
\newtheorem{remark}[theorem]{Remark}

\numberwithin{equation}{section}

\begin{document}
	
	\title{Square Root Computation in Finite Fields}
	
	\author{Ebru Adiguzel-Goktas}
	\address{Department of Computer Science, Abdullah Gul University, Kayseri, Turkey}
	
	\email{ebru.adiguzel@agu.edu.tr}
	
	\author{Enver Ozdemir}
	\address{Informatics Institute, Istanbul Technical University, Turkey}
	\email{ozdemiren@itu.edu.tr}
	
	%    General info
	\subjclass[2020]{Primary 11Y99, 68Q99}
	
	%\date{January 1, 2001 and, in revised form, June 22, 2001.}
	
	
	\keywords{Square Roots, Singular Curves, Elliptic Curves}
	
	\begin{abstract}
		In this paper, we present a review of three widely-used practical square root algorithms.  We then describe a unifying framework where each of these well-known algorithms can be seen as a special case of it. The framework with singular curves offers a broad perspective to compare and further improve the existing methods in addition to offering a new avenue for square root computation algorithms in finite fields.

	\end{abstract}
	
	\maketitle

	\specialsection*{INTRODUCTION}
	
	Finding square roots in finite fields has been an interest of many researchers in computational number theory\cite{CohenCom}. In a strict sense, the problem is reduced to finding square roots modulo a prime number $p$. In other words, the problem is converted to computing square roots in a finite field $\mathbb F_p$ where $p$ is a prime integer and $\mathbb F_p$ is a field with $p$ elements. The square root computation is a dominating step of several other methods in computations. For example, any sieving method for factorization of integers or the elliptic curve primality proving  method requires an enormous number of square root computations \cite{CohenCom}. Even though frequently a new method for computing square root in finite fields appears ({\cite{Muller}, \cite{Ozd}), the two oldest methods, Tonelli-Shanks and Cipolla ({\cite{Tonelli}, \cite{Cipolla}}), are still the most popular algorithms in practice in addition to a relatively recent one presented by Peralta \cite{Peralta}. Each algorithm has its own advantages and disadvantages over the others. First of all, all  three algorithms are probabilistic and in terms of probability of success for each trial, they are all the same considering the worst-case scenario. In other words, the probability of success for a single trial is around $1/2$ for each of these algorithms. However, considering a general case, Peralta's algorithm has higher probability of success. On the other hand, Tonelli-Shanks algorithm has lower computational load most of the time.   
		We are going to present a  method for square root computing which also allows us to observe aforementioned algorithms in the same context.  The observation leads to an efficient way of analyzing and comparing of algorithms. For example, despite Peralta's algorithm success rate being stated to be at least 1/2 in the original manuscript, we show that the actual probability can be improved to be at least 3/4 and this can be achieved by adding one step while implementing  Peralta's algorithm with a singular cubic.   We present that the success rate for a prime $p$ where $p-1=2^em$ and $m$ is odd depends on the exponent $e$.\\ 		
		\indent The well-known three algorithms are quite efficient in practice but the scholarly work on finding a method for computing square roots continues in the direction of finding a practical and deterministic algorithm \cite {Ozd}. In what follows,  we aim to explain these three algorithms as part of single method with singular curves to provide a wider insight that might help for further development toward a practical and deterministic square root algorithm. The singular curve explanation of these methods  might also help understanding of current methods and their efficiency.  In other words, each method is an analogue of finding special torsion points of certain curves where Tonelli-Shanks and Cipolla's algorithms seek sufficient conditions while Peralta's algorithm only seeks to fulfill the necessary condition for the same purpose.  In this respect, we first present a brief description of Tonelli-Shanks, Cipolla, and Peralta's algorithms in the next section. The third section  introduces the mathematical objects for the method and presentation of each algorithm in terms of these mathematical objects. We include performance analysis in the last section. 
		
		\section{Square Roots Algorithms}		
		
		Let $p$ be an odd prime integer and $\mathbb F_p$ be a prime field with $p$ elements. Let $a\neq  0 \in \mathbb F_p$. It is not hard to determine whether $a$ has a square root in  $\mathbb F_p$. In fact, the multiplicative group $G=(\mathbb F_p^*,\cdot)$ is cyclic and so $a=g^{\beta}$ for some positive integer $\beta$, where $g\in G$ is a generator of $G$.  $\beta$ is an even integer if and only if $a$ has a square root. In other words,  $\beta$ being odd implies that $a$ is a quadratic non-residue modulo $p$. It is not a tedious task to determine whether $\beta$ is even or not. As $g$ is a generator of $G,$ then 
		$$ g^{\frac{p-1}{2}}=-1\in G.$$ In this respect;
		$$a^{\frac{p-1}{2}}=(g^{\beta})^{\frac{p-1}{2}}=(g^{\frac{p-1}{2}})^\beta=(-1)^{\beta} \text { in } G.$$
		
		Therefore, $a^{\frac{p-1}{2}} \mod p$ is 1 if and only if $a$ is a quadratic residue modulo $p$. The formulation of quadratic residue/non-residue is given by the Legendre symbol;
		\[
		\left(\dfrac{a}{p}\right) =\begin{cases}
			1, & \text{if a is quadratic residue modulo p and}\, a \not \equiv 0\, \mod p \\
			-1, & \text{if a is a non-quadratic residue modulo p}\\
			0, & a\equiv 0 \mod p
		\end{cases}\]
		%That is, the Legendre symbol shows either the element $a$ of $\mathbb{F}_p$ has a square or not.
		Once decided that $a$ has a square root, the next task is to find it. 
		
		Let's assume for a moment that the prime integer $p\equiv 3 \mod 4$. We are also assuming $a$ has a square in $\mathbb F_p$, i.e., $\left(\dfrac{a}{p}\right)=1$. As
		$$\left(a^{(p+1)/4}\right)^2=a^{(p+1)/2}\equiv a^{(p-1)/2}\cdot a\equiv 1\cdot a \mod p,$$
		$x=a^{(p+1)/4}$ gives a square root of $a\in \mathbb F_p$.\\

		We do not have a generic square root formula for all other primes except for   $p\equiv 5 \mod 8 $. In fact, $a^{(p-1)/2}\equiv 1 \mod p $ implies either $$\text{  } a^{(p-1)/4}\equiv 1 \mod p \text{ or } a^{(p-1)/4}\equiv -1 \mod p.$$
		\begin{itemize}
			\item If $a^{(p-1)/4}\equiv -1 \mod p$, $x \equiv 2a(4a)^{(p-5)/8} \mod p$ gives a solution for\\ $x^2\equiv a  \mod p$  as $2^{\frac{p-1}{2}}=\left(\dfrac{2}{p}\right)= (-1)^{\frac{p^2-1}{8}}=-1$.
			\item Otherwise, $  x\equiv a^{(p+3)/8} \mod p$ gives a square root of $a$.
		\end{itemize}
		
		From now on, we are assuming all primes $p\equiv 1 \mod 8$.  The first method for square root computation that we are going to present briefly is due to Tonelli. The method is later improved by Shanks.
		\subsection{Tonelli-Shanks Algorithm} 	The multiplicative group $G=(\mathbb F^*_p,\cdot)$ has order exactly $p-1$ so for any element $s\in G$ we have $$s^{p-1}\equiv 1 \mod p.$$ 	Lets write  $p-1=2^e m$  for some $e\geq3$ and odd integer $m$.   
		Then the group $G$ has a $2$-Sylow subgroup $H$ of order  $2^e$.  Let $z$ be a generator of  $H$ and  $b$ be equal to $a^m$.} Then $b$  has order $2^j$ for some $j\ge 0$ since  $b=a^m$ lies in  $H$. As $z$ is a generator of $H$, there must be an integer $r\geq0$ such that $b\equiv z^r \mod p$. Since a is a quadratic residue, i.e., $a^\frac{p-1}{2} \equiv 1 \mod p,$ so $b$ is also a quadratic residue and then $r$ must be even integer.
	
	The fact that $r$ is even implies  $x=a^{(m+1)/2}z^{-r/2}$  a  desired square root of $a$. Indeed,
	$$x^2\equiv a^{m+1}z^{-r}\equiv a^m\cdot a \cdot b^{-1} \equiv b \cdot a\cdot b^{-1} \equiv  a \mod p$$
	Now, the problem is how to find a generator $z$ of the $2$-Sylow subgroup $H$ of $G.$ The best way of finding the generator of $H$ is to choose a random  integer $n$ such that $\left( \dfrac{n} {p}\right)=-1.$  Then $z= n^m$ is a generator of $H$ since  $z^{2^{e-1}}\equiv n^{(m2^e)/2} \equiv n^{(p-1)/2}\equiv -1 \mod p$ i.e., $z$ has order exactly $2^e$ in $H$. At the current state of the art, the only practical way to find such $n$ is the trial  and error method. Assuming GRH, $n$ is smaller than $2\log^2 p$ \cite{Bach}. Let $n$ be an element of $G$, as $G$ is cyclic and $g$ is a generator of it, there is a positive integer $\alpha$ such that $n=g^{\alpha}$. The random number $n$ gives a generator of $H$ if and only if the positive integer $\alpha$ is odd. In other words, for random $n \in G$, the probability that $n^m$ is a generator of $H$ is around $1/2$. This part of Tonelli-Shanks algorithm is probabilistic and the success rate for a single trial is almost 1/2. The resulting probabilistic algorithm for a square root computation performs the following steps:  
	\begin{algorithm}[H]\label{Tonelli}
		\caption{\textbf{: Tonelli-Shanks Algorithm}}
		\begin{algorithmic}
			\STATE\textit{Input}: a quadratic residue $a$ modulo $p$ where $p$ is an odd prime such that\\ $p-1=2^e m.$\\
			\STATE\textit{Output:} $\sqrt a \mod p$\\
			\STATE(1) Choose numbers $n$ at random until $\left( \dfrac{n} {p}\right)=-1$
			\STATE (2) Set $z=n^m \mod p$ and $b\equiv a^m \mod p.$ \\
			
			\STATE (4) Find the smallest integer $r\geq 0$ such that $b\equiv z^{r} \equiv n^{mr} \mod p.$\\
			\STATE (5) Set $x\equiv a^{(m+1)/2}z^{-r/2}\equiv \sqrt {a} \mod p.$
		\end{algorithmic}
	\end{algorithm}
	
	We should also note that finding  $r$ might be a tedious task if $e$ is a large number, say more than 99. In such a case, the usage of Tonelli-Shanks might not be suitable for practical purposes. However, Cipolla's algorithm stands as a little more practical alternative if $e$ is large.
	\subsection{Cipolla's Algorithm}
	Cipolla's Algorithm \cite{Cipolla} is described as a square root method while using  the extension  field $\mathbb {F}_{p^2}$ of $\mathbb F_p$. Even though the probability of success is the same as the Tonelli-Shanks algorithm, computing in extension brings an extra burden for this method. As mentioned above, step 4 of Algorithm 1 handles the discrete logarithm problem in the 2-Sylow subgroup $H$ of $G$. Once the size of $H$ gets large, Cipolla's method becomes more practical. In what follows, we give brief details of Cipolla's method.\\
	\indent Let $t$ be an integer with $0\leq t \leq p-1$ such that $u=t^2-a$ is  quadratic non-residue  modulo $p$. The polynomial $x^2-u$ is an irreducible polynomial in $\mathbb{F}_p$ and $\mathbb{F}_p[x]/(x^2-u)$ is isomorphic to the extension field $\mathbb{F}_{p^2}.$ Let $w=\sqrt{u}$ in $\mathbb F_{p^2}.$ As $w$ can not be an element of $\mathbb{F}_p$, one can  define the extension field as;  $$\mathbb{F}_p[\sqrt{u}]=\mathbb{F}_p[w]=\{x+yw : x,y \in \mathbb {F}_p\}.$$

	\noindent For every element $x+yw \in \mathbb{F}_p[w]$, we have $(x+yw)^p=x-yw$, since 	$$w^{p-1}\equiv (w^2)^{(p-1)/2}\equiv u^{(p-1)/2}\equiv -1 \mod p.$$
	%All binomial coefficients $\left( \dfrac{p} {j}\right)$ for $1 \leq j\leq p-1$ are divisible by $p$ and $\forall x,y \in \mathbb{F}_p$ we have $x^p=x$, $y^p=y.$ So,
	%$$(x+yw)^p=x^p+y^pw^p=x+yw^p=x-yw.$$
	\begin{theorem}\label{Cipolla}
		The element $b=(t+w)^{p+1/2}$ lies in $\mathbb{F}_p$ and it is a square root of $a$. 
	\end{theorem}
	
	\begin{proof}

		$$b^2=(t+w)^{p+1}=(t+w)(t+w)^p=(t+w)(t-w)=t^2-w^2=t^2-(t^2-a^2)=a.$$
	\end{proof}
	
	The first aim of  Cipolla's method is to find a primitive element $u$ in $\mathbb F_{p^2}$ such that $u=t^2-a$. For our case, $u$ is a primitive element if and only if $t^2-a$ is a quadratic non-residue modulo $p$. As in the case of  Tonelli-Shanks algorithm, there is no deterministic method for finding such a primitive element. The desired primitive element is found via the trial  and error method. Again, for a random $t$, the probability that $u=t^2-a$ gives a primitive element is around $1/2$. The implementation of Cipolla's algorithm follows the following steps: 
	\begin{algorithm}
		\caption{\textbf{: Cipolla's Algorithm}}
		\begin{algorithmic}	
			\STATE \textit{Input:} an odd prime $p$ and a quadratic residue $a$ modulo $p$.\\
			\STATE \textit{Output:} $\sqrt a \mod p$\\
			\STATE (1) Find an integer $t$ with $0 \leq t \leq p-1$ such that $u=t^2-a$ is a quadratic non-residue $\mod p.$\\
			\STATE (2) Return $(t+\sqrt{u})^{(p+1)/2}.$\\
		\end{algorithmic}	
	\end{algorithm}
	
	As mentioned above the last step of the algorithm requires more multiplication operations than the Tonelli-Shanks method as $w$ does not lie in $\mathbb F_p.$ It lies strictly in the extension field $\mathbb F_{p^2}$. The running time of the algorithm is bounded by $O(M\log^2 p)$ while Tonelli-Shanks is $O(M(\log p+e^2))$ where $M$ stands for complexity of multiplication in $\mathbb F_p$.    Moreover, the complexity of the Tonelli-Shanks approach can be improved by using a faster algorithm for the discrete logarithm  which  yields a complexity of $O(M(\log p+e \log e /\log \log e))$\cite{Sutherland}. In practice, this approach will almost always outperform Cipolla unless $e$ is very large.
	\subsection{Peralta's Algorithm} A relatively new and practical method by Peralta \cite{Peralta} also exploits the group  $(\mathbb Z_p^*,\cdot)$. Unlike Tonelli-Shanks and Cipolla's algorithms, the probability of success does not stand  same all the time. If $p-1=2^em$ with $m$ odd, the probability of success depends on $e$ and it is at least 1/2.
	Consider a square $a$ in the group $(\mathbb Z_p^*,\cdot)$. Let $f(y)=y^2-a$ be a polynomial over $\mathbb {F}_p.$ Note that $f(y)$ is reducible. The ring $R=\mathbb{Z}_p[y]/(f(y))$ is a factor ring and every nonzero element in $$R=\mathbb{Z}_p[y]/(f(y))=\{b_1 y+b_0+(f(y))\, | b_0, b_1 \in \mathbb{Z}_p\}$$  is either a unit or a zero-divisor.  Since the set of units $R^*$ in $R$ form a finite abelian group under multiplication modulo $f(y)$ and every finite abelian group is a direct products of cyclic groups of prime order, $R^*$ is isomorphic to $\mathbb{Z^{*}}_p\oplus\mathbb{Z^{*}}_p.$ Hence, $$e^{p-1}\equiv 1 \mod p\,\, \text{for every e} \in R^*.$$ 
	The ring $R$ can also be considered as the ring $\mathbb Z_p[\sqrt{a}]$ and the norm of any unit in this ring should be $1.$ In other words, if $r+s\sqrt{a}$ is a multiplicative unit then $N(r+s\sqrt{a})=(r+s\sqrt{a})(r-s\sqrt{a})=r^2-s^2a=1$  which also means that $r^2-s^2a\neq 0$ where  $N$ is the norm function.
	Therefore, $R^*$ consists of elements of the form $r+s\sqrt{a}$ such that $r^2 \not \equiv s^2a \mod p.$ If $r \in \mathbb{Z}^{*}_p$ and $r^2 \equiv a \mod p$, then $r$ would be our desired result. From now on, we will consider the case $r^2 \not \equiv a \mod p$ and we focus on the elements of $R^*$ of the form      $r+\sqrt{a}$. Consider $$(r+\sqrt{a})^{(p-1)/2}=u+v\sqrt{a}.$$ then we have  $$(r+\sqrt{a})^{(p-1)}\equiv 1\equiv u^2+2uv\sqrt{a}+v^2a \mod p.$$ It implies that $2uv=0$, i.e., either $u$ or $v$ is $0.$ In the case  $u=0$ we have, $$(r+\sqrt{a})^{(p-1)/2}=0+v\sqrt{a} \implies (r+\sqrt{a})^{p-1}\equiv 1\equiv v^2a \mod p.$$
	
	Then $v^{-1}$ is our desired square root of $a$. Hence, if one can find a random $r \in \mathbb{Z}_p^*$ such that $r^2 \not \equiv a \mod p$ and $(r+\sqrt{a})^{(p-1)/2}=0+v\sqrt{a}$, then the solution of the square root problem will be $v^{-1}.$ The algorithm based on this above observation follows the steps: 
	\begin{algorithm}[H]
		\caption{\textbf{: Peralta's Algorithm I}}
		\begin{algorithmic}
			\STATE \textit{Input:} a quadratic residue $a$ modulo $p$\\
			\STATE \textit{Output:} $x\equiv \sqrt a \mod p.$\\
			\STATE (1) Choose $r\in \mathbb{Z}^{*}_p$ at random such that $r^2 \not \equiv a \mod p$, otherwise output is $r.$\\
			\STATE (2) Compute $(r+\sqrt{a})^{(p-1)/2}=u+v\sqrt{a}.$\\
			\STATE (3) If $u=0$, output $x \equiv v^{-1} \mod p$ else go to $(1)$.
		\end{algorithmic}
	\end{algorithm}
	A slightly modified version of the above algorithm:
	\begin{algorithm}[H]
		\caption{\textbf{: Peralta's Algorithm II}}
		\begin{algorithmic}
			\STATE\textit{Input:} a quadratic residue $a$ modulo $p.$\\
			\STATE\textit{Output:} $x \equiv \sqrt a \mod p.$\\
			\STATE	(1) Choose $r$ at random $\in \mathbb{Z}^{*}_p.$\\
			\STATE	(2) If $r^2 \equiv -a \mod p$, choose a new $r.$\\
			\STATE	(3) Compute $(r+\sqrt{-a})^m=u+v \sqrt{-a},$ where $p-1=2^e m$ and $m$ is an odd integer.\\
			\STATE	(4) If either $u$ or $v$ is $0$, choose a new $r.$\\
			\STATE	(5) Compute $(u+v\sqrt{-a})^{2^i}$ for some $i=1,2,\cdots ,e $ until $(u+v\sqrt{-a})^{2^i}=0+w\sqrt{-a}$ for some $w.$\\
			\STATE	(6) Let $(u+v\sqrt{-a})^{2^{i-1}}=k+l\sqrt{-a}.$ Then $k^2-l^2a \equiv 0 \mod p$ and output $k/l.$
		\end{algorithmic}
	\end{algorithm}
	The probability of success depends on the number $e$. In fact, for a random $r$, we have $1-1/2^{e-1}$ chance to end up with a square root of $a$ (See \cite[Theorem 4]{Peralta}). On the other hand, as  mentioned above, when we look at Peralta's algorithm through singular curves, we would notice that the actual probability of success is $1-1/2^e$. Peralta's algorithm serves for all integers but as we do  have a generic square root formula for $p\equiv 3\mod 4$, the probability of success was defined to be 1/2 by Peralta considering the worst-case scenario ($e=2$)  event though we see in a moment that it is 3/4.
	
	\section{Geometric Analogue of  Algorithms}
	In this section, we present analogues of all three algorithms for square root computing \cite{Cipolla},\cite{Peralta},\cite{Tonelli} using  the Jacobian group of  singular cubics. In the first part, we describe the group structure of a singular  curve's Jacobian group and then show that the square root is embedded in certain points in this group \cite{Ozd11}.
	\subsection{Square Root Algorithm With Curves:}
	Let $E$ be an elliptic curve over a finite field $\mathbb F_p$ with $p$ elements where $p$ is an odd prime integer. In particular, the curve $E$ can be defined as  \begin{equation}\label{eq}
		y^2=x^3+Ax+B \text{ over the field } \mathbb F_p. 
	\end{equation}
	\subsection{Computing In An Elliptic Curve Group}
	\subsubsection{Affine Coordinates:} In case the affine coordinates is employed for representation of elements in the elliptic curve $E:y^2=x^3+Ax+B$ group, an addition operation can be performed as follows:
	
	\textbf{Addition:}	\\
	Let $P_1=(x_1,y_1), P_2=(x_2,y_2)$ be points on $E$ such that $P+Q=P_3=(x_3,y_3).$
	\begin{itemize}
		\item If $P_1\neq \pm P_2 $ then $x_3=k^2-x_1-x_2,$ $y_3=k(x_1-x_3)-y_1,$ $k=\frac{y_1-y_2}{x_1-x_2}.$
		\item If $x_1=x_2$ but $y_1 \neq y_2,$ then $P_1+P_2=\infty.$
	\end{itemize}
	
	\textbf{Doubling:}\\
	Let $2P_1=(x_3,y_3).$ Then \begin{equation}
		x_3=k^2-2x_1, y_3=k(x_1-x_3)-y_1, k=\frac{3x_1^2+A}{2y_1}.
	\end{equation}
	An addition and doubling operation require $(I+2M+S)$ and $(I+2M+2S)$ calculations, where $I$ is inversion, $M$ is multiplication and $S$ is squaring in finite field,  respectively.[\cite{Washington}]
	
	Since the calculation of inversion is costly, representing elements with projective coordinates in $E$ might be more suitable in practice.\\
	\subsubsection{Projective Coordinates:}
	For representing elements in projective coordinates, the equation of $E$ is written as
	$$Y^2Z=X^3+AXZ^2+BZ^3.$$
	The point $P_1=(X_1:Y_1:Z_1)$ on E corresponds to the affine point $(X_1/Z_1,Y_1/Z_1)$ when $Z_1 \neq 0,$ otherwise it represents the point at infinity $P_{\infty}=(0:1:0).$
	Also, $-P_1=(X_1:-Y_1:Z_1).$\\
	\indent \textbf{Addition}\\
	Let $P_1=(X_1:Y_1:Z_1), P_2=(X_2:Y_2:Z_2)$ on and $P_3=P_1+P_2=(X_3:Y_3:Z_3).$
	\begin{itemize}
		\item If $P_1 \neq P_2,$ then set
		$$K=Y_2Z_1-Y_1Z_2, L=X_2Z_1-X_1Z_2,M=K^2Z_1Z_2-L^3-2L^2X_1Z_2$$
		so that
		$$X_3=LM, Y_3=K(L^2X_1Z_2-M)-L^3Y_1Z_2, Z_3=L^3Z_1Z_2.$$
		\item If $P_1=-P_2,$ then $P_1+P_2=\infty.$
	\end{itemize}
	\indent \textbf{Doubling:}\\
	Let $2[P_1]=(X_3:Y_3:Z_3)$ then set
	$$K=4X_1Y_1^2, L=3X_1^2+AZ_1^4$$
	and
	$$X_3=-2K+L^2,Y_3=-8Y_1^4+L(K-X_3), Z_3=2Y_1Z_1.$$
	
	The complexities for an addition and doubling are $12M+4S, 4M+6S,$ respectively. No inversion is needed.

	\indent The idea of using singular curves for square root finding is based on the following observation. Let $P_1, P_2, P_3$ be points on the elliptic curve such that $P_1=(x_1,0)$, $P_2=(x_2,0)$ and $P_3=(x_3,0)$. In other words, $x_1,x_2,x_3$ are the $x$-coordinates where the graph of $E$ is crossed with the $x$-axis and this means that the tangent line at the point $P_i$ is vertical for $i=1,2,3$. The group operation in the elliptic curve group $E(\mathbb F_p)$ indicates that $P_i$ is equal to its inverse as the tangent line is vertical at these points. That implies $2P_i=\infty$ where $\infty$ represents the identity of the group $E(\mathbb F_p)$.   Since the right side of Equation \ref{eq}  can have at most three distinct roots, there are at most three  $2$-torsion points in the curve $E$'s group. The observation basically says reaching a 2-torsion on the curve gives a root of $x^3+Ax+B$.\\
	\indent In particular, define $E: y^2=(x-a)(x^2-a)$ on $\mathbb F_p$ where $a$ is a quadratic residue modulo $p$.  The equation $(x-a)(x^2-a)$ has three distinct roots in $\mathbb F_p$.   Hence, the points $P_1=(a,0)$, $P_2=(\sqrt{a},0)$ and $P_3=(-\sqrt{a},0)$ are the only $2$-torsion points of $E.$ In this regard, two $2$-torsion points on $E$ give a square root of $a.$ Since the points $\{P_1,P_2,P_3,\infty\}$ forms a subgroup of order $4$, the order of the group on the elliptic curve should be divisible by 4. Therefore, we can say that $\# E(\mathbb{F}_p)=2^em$ where $m$ is an odd integer and $e\ge 2$. We should note here that finding the number of points on the curve can take at most polynomial time \cite{CohenFrey} so we can assume that the order of $E$ is computed in an efficient way. Consider a random point $Q$ on $E$. Once $mQ$ is not the identity in $E(\mathbb F_p)$, $2^jmQ$ gives a two torsion for some $0 \le j <e$. %In some sense, the idea is a geometric analogue of Tonelli-Shanks square root algorithm. 
	One can define a square root algorithm based on this observation. However, the problematic part is locating a random point on the curve $E$. Even though, one can easily say $Q_1=(0,a)$ and $Q_2=(1,1-a)$ are points on the curve, the case $mQ_1=mQ_2=\infty$  requires locating another random point to be able to continue with the algorithm. Therefore, using a plain form of the method might stop most of the time as finding a point on the curve also requires square root computation. In order to avoid this issue, singular curves can be  used.
	\subsection{Singular Curves}
	\indent The curves for our purposes are defined in a special format that resembles the well-known elliptic curves.  Let $K$ be a field with a characteristic different from 2, the curves that we work with are defined by an equation of the form $y^2=xf^2(x)$ over $K$ where $f(x)$ is a monic square- free polynomial. Even though, the normalization of such curves returns the basic geometric tool of the projective line, the Jacobian group of these curves has the potential to be utilized for problems in computational number theory \cite{Ozd, Ozd1, Ozd2}. %The attached abstract group to such singular curves is called the Generalized Jacobian group but for the sake of simplicity, we keep using only Jacobian group for singular curves too. 
	The employed singular curves are defined in a special format.  A singular cubic curve that is utilized for square root computation is defined by $E: y^2=x(x+a)^2$ over the field $\mathbb F_p$. 	The non-singular points on $E(\mathbb F_p)$ form an abelian group and the group operation is performed as described above for the elliptic curves \cite[Section 2.10]{Washington}. Note that if a change of coordinates $x=x_1-\frac{2a}{3}$ on an equation $y^2=x(x+a)^2$ is applied then the equation becomes $y^2=x_1^3-\frac{a^2}{3}x_1-\frac{2a^3}{27}$. Therefore, we can  apply above group operation formulas for $E:y^2=x_1^3+Ax+B$  where $A=-\frac{a^2}{3}$ and $B=-\frac{2a^3}{27}$ in $\mathbb F_p$.  For the sake of simplicity, we use affine coordinates for the proofs below, however the projective coordinates  give better performance in practical applications. Since the point  $(-a,0)$ is a singular point, the point $(0,0)$ is the only  $2$-torsion point on this curve. For any point $P\neq (0,0)$ on $E,$ if $mP$ is neither  $\infty$ nor $(0,0)$ then the order of $P$ must be divisible by $4$ and $2^imP$ is a $4$-torsion point of $E$ for some $0<i<e$.  We show in a moment that the $4$-torsion points on this curve give the square root of $a$. In what follows, we show that computing square root via finding a 4-torsion is actually what other algorithms implicitly do.
	\subsubsection{Tonelli-Shanks and Peralta's Algorithms with Singular Curves}
	
	Let $a$ be a quadratic residue modulo $p$. We define a special singular cubic curve $E$ with the equation $y^2=x^2(x-a)$. As we mentioned above the normalization of $E$ at the singular point $(a,0)$ returns the projective line. This can also be seen via the following re-formalization 
	\begin{equation}\label{Per}
		\left(\frac{y}{x}\right)^2=x-a	
	\end{equation}
	
	This representation strikes the observation that leads to Peralta's algorithm. In other words, the coordinate ring of $E$ and the ring in Peralta's algorithm match. Even though, the change of coordinates of any singular cubic curve  returns to the same coordinate ring, the concealed singular cubic in Peralta's algorithm is revealed via the resulting algorithm's behavior. We are now going to show that Peralta's algorithm is nothing but the method of finding 4-torsion points in the Jacobian of $E$ where $E:y^2=x(x+a)^2$.  Note that the change of coordinates $(x=x_1+a)$ on $E$ states that $E$ can also be defined by $y^2=x^2(x-a)$.  
	
	\begin{remark}
		If an element $x \in \mathbb F_p$ is not a quadratic residue then $y^2=x(x+a)^2$ has no solution modulo $p$. Because of this, for all points $(x,y)$ on the given curve the first coordinate must be a quadratic residue. This observation allows for the parametrization $(t^2,t(t^2+a))$.
	\end{remark}
	\begin{theorem} \label{Per1}
		Let $a$ be a square in $\mathbb{F}_p,$ $ E:y^2=x(x+a)^2$ be a singular curve over $\mathbb{F}_p$. Then  the points on $E$ of order 4 give the square root of $a.$
	\end{theorem} 
	\begin{proof}Let $P=(x_1,y_1)=(t^2,t(t^2+a))$ be a point on $E.$ We are going to find $2P=(x_3,y_3)$ in terms of $t$. The tangent line $ y=k(x-x_1)+y_1$ intersects $E: y^2=x(x+a)^2=x^3+2ax^2+xa^2$ and we have, $$(k(x-x_1)+y_1)^2=x(x+a)^2. \implies 0=x^3-(k^2-2a)x^2+\cdots$$
		Since the sum of three roots of the given cubic equation is $k^2-2a,$ we have
		$$x_1+x_1+x_3=k^2-2a \implies x_3=k^2-2a-2x_1, y_3=k(x_1-x_3)-y_1. $$
		Since
		$$2y\frac{dy}{dx}=(x+a)(3x+a)$$ then we have
		$$ k=\frac{dy}{dx}=\frac{(x+a)(3x+a)}{2y}$$
		This implies
		$$x_3=k^2-2a-2x_1=\frac{(t^2-a)^2}{(2t)^2},  y_3=\frac{(t^2+a)^2(t^2-a)}{(2t)^3}.$$
		
		\noindent	It follows that; \begin{equation}
			2P=\left(\frac{(a-t^2)^2}{4t^2}, \frac{(t^2+a)^2 (t^2-a)}{8t^3}\right).
		\end{equation}
		
		\noindent  	If the order of the point $P$ on $E$ is 4 , then  $4P=\infty$ i.e., $(x_3,y_3)=2P=-2P=(x_3,-y_3)$ which implies  $y_3=0$ i.e, $(t^2+a)^2(t^2-a)=0.$ It follows that $t=\pm \sqrt{a}$ or $t=\pm \sqrt{-a}.$  $(-a,0)$ is a singular point, it is not in group and so  $t$ can not be $\pm \sqrt{-a}.$ Therefore, the points that have order 4 are $(a,2a\sqrt{a}), (a,-2a\sqrt{a})$ as $t$ can only be $\sqrt{a}$ or $-\sqrt{a}$. 
	\end{proof}
	
	Theorem \ref{Per1} suggests a way of finding square root of $a$ in $\mathbb F_p$. 
	Let $R$ be a random point on the curve $E$. The singular cubic curve $E$ has attached abelian group consisting of non-singular points in addition to the point at infinity \cite[Theorem 2.31]{Washington}. Counting the points on $E$ modulo p depends on the following observation: Every quadratic residue $x$ other than $-a$ generates two points. That means there exists $p-3$ remaining points. Adding the points $(0,0)$ and $\infty$  yields $p-1$ points. The group $E(\mathbb F_p)$ has order exactly $p-1$. Let
	$p-1=2^em$ for some odd integer  $m$  and a positive integer  $e \geq 2$. Let $R=(\ell^2,\ell(\ell^2+a))$ be a random point on $E(\mathbb F_p)$ for some $\ell\neq 0\in \mathbb F_p$. If  $T=mR$ and 2T are not the identity then $2^iT$ is definitely a 4-torsion point in $E(\mathbb F_p)$ for some $0 \leq i\leq e-1$. The above discussion leads to us the following algorithm which is implicitly the method presented by Peralta \cite{Peralta}. 
	\begin{algorithm}[H] \label{SingPer}
		\caption{}
		\begin{algorithmic}
			\STATE \textit{Input:} An odd prime number $p$ and a quadratic residue $a$ modulo $p.$\\ $E:y^2=x(x+a)^2$ over the field $\mathbb F_p.$ Assume $\#E(\mathbb F_p)=2^em$ such that $(m,2)=1.$\\
			\STATE \textit{Output:}$\sqrt{a} \mod p.$\\
			
			\STATE (1) Choose a point $R=(t^2,t(t^2+a))$ on $E$ such that $Q=mR\neq \infty$ and $Q=mR \neq (0,0).$\\
			\STATE (2) $2^iQ=(z,w)$ must be a 4-torsion (i.e $z=a$) for some $i=0,1,\cdots e-1.$\\
			\STATE(3) Compute $w/2a$ which gives $\sqrt a \mod p.$\\
		\end{algorithmic}
	\end{algorithm}
	
	%Note that, if $mP=(0,0)$ in the first step, then the point $Q=(\frac{m+1}{2})P$ is a 4-torsion point and gives the square root of $a.$
	\begin{theorem}
		Let $E$ be as above and $R$ be a random point on it. The probability of reaching a point of order 4 via computing $2^imR$  for some $0\leq i < e-1$ is $1-\dfrac{1}{2^{e-1}}$.
	\end{theorem}
	
	\begin{proof}
		Since the attached abelian group, $Jac(E)$, to $E$ is cyclic group of order $p-1,$ where $p-1=2^e m,$ $Jac(E)\cong \mathbb{Z}_{2^e} \oplus \mathbb{Z}_m$. The probability that a random point's order is not divisible by 4 is $2/2^e=1/2^{e-1}$ and this completes the proof.	
	\end{proof}

	The same probability of success reveals that the curve is  the actual singular cubic curve that Peralta's algorithm employed.  The following lemma indicates an additional step can improve the success rate of Peralta's algorithm.
	\begin{lemma}\label{lemma1}
		Let $E$ be as above over the field $\mathbb F_p$ and $R=(\ell^2,\ell(\ell^2+a))$ be a random point on $E(\mathbb F_p)$. If $mR=(0,0)$ then the point $Q=(\frac{m+1}{2})R$ has order divisible by 4 where again $p-1=m2^e$ for some positive integer  $e$ and odd integer $m$.
		
	\end{lemma}
	
	\begin{proof}
		The point $Q=(\frac{m+1}{2})R$ is  in the group $E(\mathbb F_p)$ where the order of $E(\mathbb F_p)=m2^e$ with odd integer $m$ and positive integer $e$. It is given that $mR=(0,0)$ where $(0,0)$ is the only point $E(\mathbb F_p)$ of order 2. Lets look at:
		$$4Q= 2mR+2R=2(0,0)+2R=\infty+2R=2R$$
		Since we have $2mQ=mR=(0,0)$, the order of $Q$ is divisible by 4 and $mQ$ has order exactly 4. 
	\end{proof}
	
	\begin{remark}\label{remark1}
		Lemma \ref{lemma1}  implies that the success rate of Algorithm 5 can easily be improved. In other words, unless one reaches $mR= \infty$ in step 1, one can {\color{red} still} find a square root of $a$. Because,  in the case of  $mR=(0,0)$ then $Q=(\frac{m+1}{2})R$ gives a point of order divisible by 4. This observation reduces the probability of failure to $1/2^e$ and increases the success rate to $1-1/2^e$ in a single trial. Considering the worst case where $e=2$, we have at least $3/4$ chance to reach a square root in a single trial with a little modification in Algorithm 5. 
	\end{remark}	
	
	\begin{theorem}
		Let $E$ be as above and $R$ be a random point on it. The probability of $mR$ being a generator of  the Sylow-2 subgroup of $E(\mathbb F_p)$ is 1/2.
	\end{theorem}
	
	\begin{proof}
		Consider the group $\mathbb Z_p^*\cong\mathbb{Z}_{2^e} \oplus \mathbb{Z}_m$ and an element $n$ of it where $p-1=m2^e$ for some odd integer $m$. The element $n^m$ is a generator of the Sylow-2 subgroup of $\mathbb Z_p^*$ if it satisfies $n^{(p-1/2)} \equiv -1 \mod p$. Similarly, if a point $R$ in the cyclic group $E(\mathbb F_p)$ satisfies $(\frac{p-1}{2})R=(0,0)$ then $mR$ becomes a generator of the Sylow-2 subgroup of $E(\mathbb F_p)$.  As the groups $\mathbb Z_p^*$ and $E(\mathbb F_p)$ are isomorphic, the probability of reaching a generator of the Sylow-2 subgroup of $E(\mathbb F_p)$ via computing $mR$ for a random point $R$ is 1/2.    
	\end{proof}
	Once one finds a generator $mR$ of Sylow-2 subgroup of $E(\mathbb F_p)$, one can easily reach a square root by computing a 4-torsion via $2^{e-2}mR$. In other words, like  Tonelli-Shanks algorithm, finding a generator is sufficient to reach a desired square root in a deterministic way. However, Peralta's algorithm shows that this is not a necessary condition as one can reach a 4-torsion without finding a generator for  the 2-Sylow subgroup of $E(\mathbb F_p)$.

	\subsubsection{Cipolla's Algorithm with Singular Curves}
	In the same context, we now present a geometric analogue of Cipolla's algorithm.
	\begin{theorem}
		Let $a$ be a square in $\mathbb F_p,$ $E:y^2=x(x+a)^2$ be a singular curve over $\mathbb F_p$ and $P=(t^2,t(t^2+a))$ be any point on $E,$ where $t\neq 0.$ If $t^2+a$ is quadratic non-residue, then employing $P$ with Algorithm 5 definitely returns a square root of $a$.
	\end{theorem}
	
	\begin{proof}
		
		We first show that if $t^2+a$ is a quadratic non-residue, then it is impossible to have $mP = \infty.$ Suppose $mP=\infty$. Consider the point $$Q=\frac {(m+1)}{2}P.$$
		Then we have $$2Q=mP+P=P.$$ Let the $x$-coordinate of $Q$ be  $ b$ where $b\in \mathbb F_p$.  Then the x-coordinates of $2Q$ would be $\dfrac{(a-b^2)^2}{(2b)^2}$  by (2.4). So, we should have $$\dfrac{(a-b^2)^2}{(2b)^2}=t^2\,\,\text{and then}\,\, a-b^2=\pm 2bt.$$ 
		Since the discriminant of $b^2 \pm 2bt-a=0$ is $\Delta=4t^2+4a=4(t^2+a)$ and $t^2+a$ is  quadratic non-residue, it is impossible to have such  $b^2$ such that $mP=\infty$.
		
		If $mP=(0,0)$ Lemma \ref{lemma1} then Remark \ref{remark1} implies that $Q=(\frac{m+1}{2})P$ has order divisible by 4 which means  a modified version of  Algorithm 5 returns a square root.
	\end{proof}	
	The above theorem indicates that Cipolla's method also seeks a sufficient condition in Algorithm 5 similar to the Tonelli-Shanks method. Employing 4-torsion points on $E$ returns a square root algorithm where each of three practical algorithms can be explained in the same context. Towards a deterministic and polynomial-time algorithm for square root computation in the finite field, it might be beneficial to work with other torsion points. For example, designing a singular curve where 5-torsion points return square root increases the success rate. 
	
	\begin{example}
		We are going to find the square root of $a=2$ modulo $p=2017$ via a singular cubic analogue of each algorithm. We define the curve $E$ by the equation $y^2=x(x+2)^2$ over the finite field $\mathbb F_{2017}$. As $$2017-1=2016=2^5 \cdot 63$$  
		$m=63$ and $e=5$.\\
		{\bf Peralta's Algorithm:} We select a random point $P$ on $E$, say $P=(1,3)$. We compute $$Q=mP=63P=(2,90).$$	
		Since $Q$ is neither the identity nor two torsion point  $(0,0),$ the order of $P$ is divisible by 4 which means that the order of $2^iQ$ is divisible by 4 for some $i=0,\dots,3$. As $2Q=(0,0)$, the point $Q$ has order 4. Note that $$Q=(a,2a\sqrt{a})=(2,4\sqrt{2})=(2,90) \text{ which implies } \sqrt{2}=1031.$$ 
		
		{\bf \noindent Tonelli-Shanks Algorithm: } We search a point $P$ such that $$\left(\dfrac{p-1}{2}\right)P\neq \infty$$
		The first random points $P=(1,3)$ and $P=(25, 135)$ do not work. The third point $P=(289, 913)$ works as $1008P=(0,0)$. That means that $mP=Q=63P=(138,258)$ is a generator of the Sylow-2 subgroup where $$8Q=8\cdot63P=(2,1927)$$ is a point of order 4.
		Note that $(2,1927)=(2,4\sqrt{2})=(2,4\cdot 986)$ which implies $\sqrt{2}\equiv 986\equiv-1031 \mod 2017$	
		
		{\bf \noindent Cipolla Algorithm:} We first search a random point $t$ such that $t^2+a=t^2+2$ is a quadratic non-residue mod 2017. 
		We have found such $t$ on the second try where we set $t=611$ and $P=(t^2,t(t^2+2))=(176,1857)$. The order of the point $P$ must be divisible by 4. We compute $mP=(1379,1791)$ then
		$2mP=(1553, 936)$, $4mP=(96,384)$ and $8mP=(2,90)$ where we get $4\sqrt{2}=90$ and $\sqrt{2}=1031$.
	\end{example}
	
	We implemented all algorithms in the same environment where we use a C++ library, PARI/GP \cite{PARI}, for operations involving large integers. The real time tests are conducted on the computer running Linux OS with an Intel i7-11370H processor and 32 GB main memory. The discrete logarithm problem in Tonelli-Shanks algorithm is handled via a naive brute force method. The first implementation of Tonelli-Shanks seeks a quadratic non-residue randomly in the field $\mathbb F_p$ for a prime integer $p$. On the other hand, we also implement Tonelli-Shanks with the quadratic reciprocity law to find a quadratic non-residue while searching only numbers less than 100. The following table briefly summarize the test results. 
	
	\begin{center}
		\begin{table}[H]
			\begin{tabular}{| l | l | l | l |}
				\hline
				Finite Field Size (size of $p$) & 256-bit, e=4 & 512-bit,e=5 & 1024-bit, e=8 \\
				\hline
				Tonelli-Shanks & 0.753  & 1.548  & 4.792   \\
				\hline
				Tonelli-Shanks (Quadratic Reciprocity) &  0.328  & 0.642 & 2.372\\
				\hline 
				Cipolla & 0.583 &1.391  &4.484 \\
				\hline
				Peralta & 0.407  &0.720 &2.188 \\
				\hline
				Singular Cubics & 0.317  &0.992  &4.298 \\
				
				\hline
			\end{tabular}
			\caption{The tests are conducted for primes $p$ where $p-1=2^em$ and the  time (in millisecond) is average of 1000 runs for each algorithm. }
		\end{table}
	\end{center}

	\section*{Acknowledgments}
	We  thank  Andrew Sutherland and  Rewievers for their corrections and suggestions to improve the quality of paper.

	\bibliographystyle{amsplain}

\end{document}